\documentclass[10pt]{paper}

\usepackage{amssymb,amsmath,enumerate,theorem,epsfig,rotating,graphics,changebar,eepic}%t1enc
\usepackage{amsfonts}
\usepackage{a4,latexsym,parskip}
\usepackage{hyperref}
\hypersetup{pdfauthor=Matthias Sch"utt} \hypersetup{pdftitle=Diss}

\newcommand{\PP}{\mathbb{P}}

\newcommand{\Q} [1] []{\mathbb{Q}_{#1}}
\newcommand{\N} [1][] {\mathbb{N}_{#1}}

\newcommand{\F}{\mathbb{F}}
\newcommand{\Z}{\mathbb{Z}}
\newcommand{\p}{\mathfrak{p}}

\newcommand{\OO}{\mathcal{O}}

\newtheorem{Spezial-Theorem}{Theorem}[section]
\newtheorem{Spezial-Proposition}{Proposition}[section]
\theoremstyle{break} \newtheorem{Theorem}{Theorem}[section]
\newtheorem{Proposition}[Theorem]{Proposition}
\newtheorem{Lemma}[Theorem]{Lemma}

\newtheorem{Corollary}[Theorem]{Corollary}

\newtheorem{Remark}[Theorem]{Remark}

\newtheorem{Conjecture}[Theorem]{Conjecture}
\ifx\undefined\pdfpageheight
  
\else
  
\fi

\ifx\undefined\pdfpageheight
  \def\imagerotate#1#2#3#4#5#6{\put(#2){\epsfig{file=#1.eps,width=#3}}}
\else
  \def\imagerotate#1#2#3#4#5#6{\put(#5){\epsfig{file=#4,width=#6}}}
\fi

\begin{document}
\setlength{\unitlength}{1cm}

\title{Arithmetic of a singular K3 surface}

\author{Matthias Sch\"utt}

\date{\today}
%\subjclass[2000]{14J27,14J28}
\maketitle

\abstract{This paper is concerned with the arithmetic of the
elliptic K3 surface with configuration [1,1,1,12,3*]. We determine
the newforms and zeta-functions associated to $X$ and its
twists. We verify conjectures of Tate and Shioda
for the reductions of $X$ at $2$ and $3$.}

\begin{small}

\end{small}

\vspace{0.2cm}

\section{Introduction}
\label{s:intro}

This paper investigates the arithmetic of a particular singular K3
surface $X$ over $\Q$, the extremal elliptic fibration with
configuration [1,1,1,12,3*]. First of all, we determine the
corresponding weight 3 form (cf.~\cite[Ex.~1.6]{L}) explicitly.
For this, we calculate the action of Frobenius on the
transcendental lattice by counting points and applying the
Lefschetz fixed point formula. The proof is based on our previous
classification of CM-forms with rational coefficients in
\cite{S-CM}. In fact, we only have to compute one trace.

Then we compute the zeta-function of the surface. This is used to
study the reductions of $X$ modulo some primes $p$. We emphasize
that we are able to find a model with good reduction at 2. We
subsequently verify conjectures of Tate and Shioda. The
conjectures will be recalled in Section \ref{s:div-conj} and
verified in Sections \ref{s:X/3} - \ref{s:X/2}.

The final section is devoted to the twists of $X$. We show that
these produce all newforms of weight 3 with rational coefficients
and CM by $\Q(\sqrt{-3})$.

\section{The extremal elliptic K3 fibration}
\label{s:3*} \label{s:rat-ell}

There is a unique elliptic K3 surface $X$ with configuration
[1,1,1,12,3*] and with a section (No.~166 in \cite{SZ} and
\cite[Tab.~2]{S-K3}). Since this arises as cubic base change of
the extremal rational elliptic surface $Y$ with singular fibres
$I_1^*, I_4$ and $I_1$, we shall start by studying this surface.

In \cite{MP1}, an affine Weierstrass equation of this fibration
was given as \begin{eqnarray}\label{eq:1} Y':\;\;y^2 =
x^3-3\,(s-2)^2 (s^2-3)\, x+s (s-2)^3 (2 s^2-9).
\end{eqnarray}
It has discriminant
\[
\Delta=16\cdot 27\,(s-2)^7(s+2),
\]
so the singular fibres are $I_1^*$ above $2, I_1$ above $-2$ and
$I_4$ above $\infty$.

We shall look for a model of $Y$ over $\Q$ which has everywhere
good reduction. The fibre of $Y'$ at $\infty$ has non-split
multiplicative reduction, so H$_{\text{\'et}}^2(Y',\Q[\ell])$ is
ramified. Therefore we twist equation (\ref{eq:1}) over the
splitting field $\Q(\sqrt{-3})$. Performing some elementary
transformations (cf.~\cite[IV.1]{S-D}), we obtain the equation
\begin{eqnarray}\label{eq:rat-ell-2}
Y:\;\;\;y^2+s\,x\,y=x^3+2s\,x^2+s^2\,x.
\end{eqnarray}
This elliptic surface has discriminant $\Delta=s^7(s+16)$ and
everywhere good reduction. In particular, all components of
reducible fibres are defined over $\Q$. (For $I_1^*$, this can be
derived from Tate's algorithm \cite[IV, \S 9]{Si}.) Upon reducing
mod 2, we obtain the equation from \cite{Ito}. $Y$ then only
inherits the two singular fibres of types $I_4$ and $I_1^*$ with
wild ramification at the latter fibre.

We now come to the singular K3 surface $X$. Consider the cubic
base change
\[
\pi: s \mapsto s^3.
\]
Via pull-back from $Y$, this gives rise to an extremal elliptic K3
surface. The resulting Weierstrass equation reads
\begin{eqnarray}\label{eq:K3}
y^2+s^2\,x\,y=x^3+2s\,x^2+s^2\,x.
\end{eqnarray}
This Weierstrass model has a $D_7$ resp.~$A_{11}$ singularity in
the fibre above $0$ resp.~$\infty$. By $X$, we denote the minimal
desingularization. This elliptic fibration has configuration
[1,1,1,12,3*]. The reducible singular fibres $I_3^*$ and $I_{12}$
sit at 0 and $\infty$. The three fibres of type $I_1$ can be found
at the cube roots of $-16$.

By construction, the absolute Galois group Gal$(\overline\Q/\Q)$
acts trivially on the trivial lattice $V$ of $X$, generated by the
$0$-section and the fibre components. Since $X$ is extremal,
tensoring $V$ with $\Q$ gives the corresponding statement for
$NS(X)$.

The elliptic surface $Y$ has everywhere good reduction.
Furthermore, the base change $\pi$ is nowhere degenerate upon
reducing. Hence, the pull-back $X$ can only have bad reduction at
the prime divisors of the degree of $\pi$, i.e.~at 3. Modulo 3,
the Weierstrass model (\ref{eq:K3}) obtains an additional $A_2$
singularity, so the reduction is in fact bad
(cf.~Sect.~\ref{s:X/3}).

In terms of H$_{\text{\'et}}^2(X,\Q[\ell])$, the ramification is
reflected in the contribution of the transcendental lattice $T_X$,
the orthogonal complement of $NS(X)$ in H$^2(X,\Z)$. Here we
consider it as a two-dimensional $\ell$-adic Galois representation
$\rho$. The reduction properties of $X$ imply that $\rho$ is only
ramified at 3 and the respective prime $\ell$. This agrees with
the discriminant $d=d_{T_X}$ of $X$ which is 3. To see this,
recall that
\[
d_{T_X}=-d_{NS(X)},
\]
since H$^2(X,\Z)$ is unimodular. Consider the trivial lattice $V$
of $X$. Let $U$ denote the hyperbolic plane, i.e.~$\Z^2$ with
intersection form $\begin{pmatrix} 0 & 1\\1 & 0\end{pmatrix}$.
Then
\[
V=A_{11}\oplus D_7\oplus U
\]
has discriminant $d_V=-12\cdot 4=-48$. The N\'eron-Severi group
$NS(X)$ is obtained from $V$ by adding the sections. Here, $X$ is
extremal, so $MW(X)$ is finite. Hence
\[
d_{NS(X)}=\frac{d_V}{|MW(X)|^2}.
\]
Explicitly, we have $MW(X)\cong MW(Y)$, consisting of four
elements (cf.~Sect.~\ref{s:NS}). We obtain $d_{NS(X)}=-3$ and
$d_{T_X}=3$. Then the (reduced) intersection form on $T_X$ can
only be
\[
\begin{pmatrix} 2 & 1\\1 & 2\end{pmatrix}.
\]

\section{The associated newform}
\label{s:new}

By a result of Livn\'e \cite[Ex.~1.6]{L}, $X$ (or $T_X$
resp.~$\rho$) has an associated newform of weight 3 with CM by
$\Q(\sqrt{-3})$. Our aim is to determine this newform explicitly.
Let $p$ be a prime of good reduction, that is $p\neq 3$. Then
\begin{eqnarray}\label{eq:det}
\text{det } \rho (\text{Frob}_p) = \left(\frac{p}{3}\right)p^2.
\end{eqnarray}
To find the trace of $\rho$, we use the Lefschetz fixed point
formula. We have already seen that Gal$(\overline\Q/\Q)$ operates
trivially on $NS(X)$. Hence, the Lefschetz fixed point formula
gives
\[
\# X(\F_p)=1+20\,p+\text{tr } \rho(\text{Frob}_p)+p^2.
\]
Using a computer program, we calculated the following traces at
the first good primes.

\begin{center}
\begin{tabular}{|c||c|c|c|c|c|c|c|c|c|c|c|}
\hline $p$ & 2 & 5 & 7 & 11 & 13 & 17 & 19 & 23 & 29 & 31 & 37\\
\hline
tr $\rho$ & 0 &  0 & $-13$ & 0 &  $-1$  & 0 & 11 & 0 & 0 &
$-46$ &
47\\
\hline
\end{tabular}
\end{center}

By inspection, these traces coincide with the Fourier coefficients
of the newform $f=\sum_n a_nq^n$ of level 27 and weight 3 from
\cite[Tab.~1]{S-CM}. We shall now prove that this holds at every
prime:

\begin{Proposition}\label{Prop:Euler}
\[
L(T_X,s)=L(f,s).
\]
\end{Proposition}

The proof makes use of Livn\'e's modularity result for $X$. Let
$g$ denote the associated newform. Since $g$ has CM by
$\Q(\sqrt{-3})$, it is a twist of $f$ by \cite[Thm.~3.4]{S-CM}. In
our special situation, $g$ is unramified outside 3 (i.e.~it has
level $3^r$), since $X$ has good reduction elsewhere. Hence we
only have to compare the few possible newforms with such a level.

Let $\psi_f$ denote the Gr\"o\ss{}encharakter of conductor 3 which
corresponds to $f$, and analogously $\psi_g$. Then
\cite[Prop.~II.11.1]{S-D} and in particular the reasoning of section
II.11.4 therein show that there are three possibilities in total:
\begin{eqnarray}\label{eq:twist-f}
\psi_g=\psi_f\;\;\text{ or
}\;\;\psi_g=\psi_f\otimes\left(\frac{3}{\cdot}\right)_3\;\;\text{
or }\;\; \psi_g=\psi_f\otimes\left(\frac{3}{\cdot}\right)_3^2.
\end{eqnarray}
Since the third residue symbol evaluates non-trivially at the
factors of $p=7$, it suffices to compare the Fourier coefficients
(or traces) at 7. Thus $f=g$. A priori, this only guarantees that
all but finitely many traces coincide. But here, the associated
Galois representations $\rho$ and $\rho_f$ are simple, since their
traces are not even. Hence Proposition \ref{Prop:Euler} follows.

\begin{Corollary}\label{Cor:zeta}
The zeta-function of $X$ is
\[
\zeta(X,s)=\zeta(s)\;\zeta(s-1)^{20}\,L(f,s)\;\zeta(s-2).
\]
\end{Corollary}

\emph{Proof:} We verify the corollary at every local Euler factor.
On the one hand, at all primes $p\neq 3$, this follows from
Proposition \ref{Prop:Euler} and $\rho(X/\Q)=20$.

On the other hand, the reduction $X_3$ at 3 is singular. Hence,
the local Euler factor is only defined via the number of points of
this singular variety over the fields $\F_{3^r}$, since the
Lefschetz fixed point formula is not available.

The idea is that, with respect to the number of points, $X_3$
looks like $\PP^2$ blown up in 19 rational points. Let us explain
what we mean by this: We compare $X_3$ to $Y_3$, the reduction of
$Y$ modulo 3. The base change $\pi$ is purely inseparable mod 3.
Hence, over any finite field $\F_{3^r}$, the smooth fibres of
$X_3$ have the same number of points as the corresponding fibres
of $Y_3$. On the other hand, $X_3$ has ten additional $\PP^1$s in
the singular fibres at 0 and $\infty$. These are all defined over
$\F_3$. Thus
\[
\#X_3(\F_{3^r})=\#Y_3(\F_{3^r})+10\cdot 3^r.
\]
Recall that $Y_3$ is $\PP^2$, blown up in nine points. These
points are all rational over $\F_3$, since the absolute Galois
group operates trivially. We deduce that
\[
\#X_3(\F_{3^r})=\#(\PP^2(19))(\F_{3^r}).
\]
This gives the local Euler factor
\[
\zeta_3(X,T)=\dfrac{1}{(1-T)(1-3\, T)^{20}(1-3^2\, T)}.
\]
The level 27 implies $L_3(f,s)=1$ by classical theory. This can
also be read off from \cite[Tab.~2]{S-CM} and the nebentypus
$\chi_{-3}=(\frac \cdot{3})$ of $f$. This completes the proof of
Corollary \ref{Cor:zeta}.

\section{The conjectures for the reductions}
\label{s:div-conj}

In this section, we shall discuss conjectures of Tate, Shioda and
Artin for smooth projective surfaces over finite fields. In
particular, these conjectures apply to (supersingular) reductions
of varieties defined over number fields.

Corollary \ref{Cor:zeta} will be very useful: If $X$ has good
reduction at $p$, then the corollary gives the local
$\zeta$-function of the smooth variety $X/\F_p$. Explicitly, let
$p\neq 3$. We obtain
\[
P_2(X/\F_p, T)=\text{det}(1-\text{Frob}_p\,T;
\text{H}_{\text{\'et}}^2(X/\overline{\F}_p,\Q[\ell]))=(1-p\,T)^{20}(1-a_p\,
T + \chi_{-3}(p)\, p^2\, T^2).
\]
This is exactly where the Tate conjecture enters. To formulate it,
consider a finite field $k$ and a smooth projective variety $Z/k$.
Define the Picard number of $Z$ over $k$
\[
\rho(Z/k)= \text{rk } NS(Z)^{\text{Gal}(\overline k/k)}.
\]
We employ the convention $\rho(Z)=\rho(Z/\overline k)$ when the
field of definition of $Z$ is understood.

\begin{Conjecture}[Tate {\cite[(C)]{T1},
\cite{T2}}]\label{Conj:Tate} Let $q=p^r$ and $Z/\F_q$ a smooth
projective variety. Denote the order of the zero of $P_2(Z/\F_q,
T)$ at $T=\frac{1}{q}$ by $u$. Then $u=\rho(Z/\F_q)$.
\end{Conjecture}

The Tate conjecture is known for elliptic K3 surfaces with a
section in characteristic $p>3$ \cite[Thm.~(5.6)]{T2}.

We can consider the Weierstrass model (\ref{eq:K3}) over any
$\F_q, q=p^r$. Denote the minimal resolution of the $A_n$ and
$D_m$ singularities by $X/\F_q$. This surface coincides with the
reduction $X_p$ of $X$ if and only if $p$ is a good prime
(i.e.~$p\neq 3$). On the other hand, $X_3$ contains an $A_2$
singularity, so it is not smooth. The desingularization $X/\F_3$
will be sketched in the next section.

\begin{Proposition}\label{Prop:Tate}
Let $p$ be a prime. Consider $X/\F_p$. Then
\[
\rho(X/\F_p)=\begin{cases}
20, & \text{if }\; p\equiv 1\mod 3,\\
21, & \text{if } \;p\equiv 0, 2\mod 3.
                \end{cases}
\]
\end{Proposition}

If $p>3$, Proposition \ref{Prop:Tate} follows from the (known)
Tate conjecture. The proofs for $p=2$ and 3 will be given in the
following three sections. We will also verify a conjecture of
Shioda which is introduced in the following.

Let $L$ be a number field and $Z$ a singular K3 surface over $L$.
If $\p$ is a prime of $L$, denote the residue field of $L$ at $\p$
by $L_\p$. We call $\p$ supersingular if and only if $Z/L_\p$ is
supersingular (i.e.~$\rho(Z/\overline L_\p)=22$).

Shioda's conjecture concerns the surface $Z/L_\p$ at a
supersingular prime $\p$. We can compare two lattices of rank two:
On the one hand, we have the transcendental lattice $T_Z$ of
$Z/L$. On the other hand we can use the natural embedding
\[
NS(Z/\overline L)\subseteq NS(Z/\overline L_\p)
\]
to define the orthogonal complement
\[
T_\p=NS(Z/\overline L)^\bot\subset NS(Z/\overline L_\p).
\]

\begin{Conjecture}[Shioda {\cite[Conj.~4.1]{ShCR}}]\label{Conj:Shioda}
Let $Z$ be a singular K3 surface over a number field $L$ and $\p$
a supersingular prime. Then, the two lattices $T_Z$ and $T_\p$ are
similar.
\end{Conjecture}
In other words, the claim is that $T_\p$ is isomorphic to
$T_Z(-m)$ for some $m\in\Q[>0]$. We will verify this conjecture
for the extremal elliptic K3 fibration $X$ at the primes 2 and 3
in the next three sections. This will be achieved by finding
explicit generators of the respective $T_\p$. As a by-product,
we will thus verify the following

\begin{Theorem}[Artin {\cite[(6.8)]{A}}]\label{Thm:Artin}
Let $Z$ be a supersingular K3 surface
over a finite field $k$ of $p^r$ elements. If $r$ is odd, then
Gal$(\overline k/k)$ operates non-trivially on $NS(Z)$.
\end{Theorem}

At the time of this paper, the Brauer group $Br(Z)$, if finite, was only known to have cardinality a square or twice a square. Hence, Artin was only able to prove the above result for odd $p$. Recently, the second alternative has been ruled out in \cite{Brauer}. Hence, Artin's argumentation in \cite[\S 6]{A} also applies to $p=2$.

Note that Theorem \ref{Thm:Artin} agrees perfectly with Proposition \ref{Prop:Tate}
for our surface $X$. In fact, the known part of Tate's Conjecture
\ref{Conj:Tate} shows that at a supersingular prime $p>3$ (i.e.~$p\equiv -1\mod 3$),
\[
\rho(X/\F_{p^2})=22.
\]

\section[$X/{\F}_3$]{$\boldsymbol{X/{\F}_3}$}
\label{s:X/3}

In this section, we shall consider $X/\F_3$. This will be special
because $\pi$ is purely inseparable modulo $3$. As a consequence,
the base change $X/\F_3$ from $Y/\F_3$ via $\pi$ has only three
singular fibres. They have types $I_3^*, I_{12}$ and $I_3$. In
particular, the elliptic fibration $X/\F_3$ is extremal (cf.~the
classification of \cite{Ito}) and supersingular. By construction, it also is unirational. We emphasize that
the reduction $X_3$ is not smooth.

We shall now verify Tate's and Shioda's conjectures for $X/\F_3$.
Consider the 4-torsion sections of $X/\F_3$ which come from $X/\Q$
upon reducing (see Section \ref{s:NS} for a detailed study). These
sections meet the $O$-component of the $I_3$ fibre. Denote the
other components of the $I_3$ fibre (not meeting $O$) by
$\Theta_1, \Theta_2$. By construction, they also do not meet the
generators of the trivial lattice $V$ of $X/\Q$ (embedded into
$NS(X/\F_3)$). Hence, they are orthogonal to $NS(X/\Q)\subset
NS(X/\overline\F_3)$. We claim that this is already all of
$NS(X/\overline\F_3)$:

\begin{Lemma}\label{Lem:X/3}
Let $A_2$ denote the root lattice generated by $\Theta_1$ and $\Theta_2$. Then
\[
NS(X/\overline\F_3)\cong NS(X/\Q) \oplus A_2.
\]
\end{Lemma}

The proof of this lemma will be directly derived from the following classical

\begin{Theorem}[Artin, Rudakov-\v Safarevi\v c]\label{Thm:Artin-inv}
Let $X$ be a supersingular K3 surface over a field of characteristic $p$. Then
\[\text{discr } NS(X)=-p^{2\sigma_0} \;\;\;\text{ for some }\;\sigma_0\in\{1,\hdots,10\}.\]
Here, $\sigma_0$ is called  the \emph{Artin invariant}.
\end{Theorem}

As a consequence, in our situation, we have discr
$NS(X/\overline\F_3)=-3^{2\sigma_0} \text{ for some
}\;\sigma_0\in\{1,\hdots,10\}$. But then, the above sublattice
$NS(X/\Q) \oplus A_2$ of $NS(X/\overline\F_3)$ has rank 22 and
discriminant $-9$. Since this is the maximal possible
discriminant, we deduce the equality of the two lattices. This
proves Lemma \ref{Lem:X/3}.

\begin{Corollary}
The supersingular K3 surface $X/\F_3$ has Artin invariant
$\sigma_0=1$.
\end{Corollary}

We shall now prove Proposition \ref{Prop:Tate} for $p=3$. Since
$\rho(X/\Q)=20$, the reduction of $NS(X/\Q)$ is clearly generated
by divisors over $\F_3$. Using Lemma \ref{Lem:X/3}, we only have
to study the field of definition of the two further generators
$\Theta_1, \Theta_2$ of $NS(X/\F_3)$.

For this, consider the elliptic curve over $\F_3(s)$ given by
equation (\ref{eq:K3}). It has non-split multiplicative reduction
at $s=-1$. More precisely, the components $\Theta_1$ and
$\Theta_2$ are conjugate in $\F_3(\sqrt{-1})$. In particular,
their sum is defined over $\F_3$, while the difference has
eigenvalue $-1$ with respect to the conjugation. Hence, we deduce
the claim of Proposition \ref{Prop:Tate} that $\rho(X/\F_3)=21$.
This agrees with the Tate conjecture and Theorem \ref{Thm:Artin}.
We also see that $\rho(X/\F_9)=22$.

Finally, we come to Shioda's Conjecture. By Lemma \ref{Lem:X/3},
we have
\[
T_3=NS(X/\Q)^\bot=A_2\subset NS(X/\overline\F_3).
\]
Since $A_2$ has intersection matrix $\begin{pmatrix} -2 & 1\\ 1 & -2\end{pmatrix}$, we deduce the validity of Conjecture \ref{Conj:Shioda} with $m=1$.

\section[$NS(X)$]{$\boldsymbol{NS(X)}$}
\label{s:NS}

To verify the conjectures for the reduction $X_2$, we need a
better knowledge of $NS(X)=NS(X/\overline\Q)$. To be precise, we
want to express the sections of the original fibration $X$ over
$\Q$ in terms of $V\otimes_\Z\Q$. This is possible since the
fibration is extremal, such that
\[
MW(X)\cong NS(X)/V
\]
is only torsion. In terms of equation (\ref{eq:K3}), the sections
of the elliptic fibration $X$ are given by:
\[
MW(X)\cong\Z/4=<P>=<(-s,0)>=\{(-s,0),(0,0),(-s,s^3),O\}.
\]
They can be derived from the sections of $Y'$, as given in
\cite{MP1}, by following them through the base and variable
changes. We will later see that $MW(X)$ gives all torsion-sections
of $X/\overline\F_2$ upon reducing.

We want to express the section $P$ as a $\Q$-divisor in
$V\otimes_\Z\Q$. This can be achieved by determining the precise
components of the singular fibres which $P$ intersects.

In the $I_{12}$ fibre above $\infty$, we number the components
cyclically $\Theta_0,\hdots,\Theta_{11}$ such that $\Theta_0$
meets $O$. The fibre of type $I_3^*$ above $0$ is sketched in
Figure \ref{Fig-3}. As usual, $C_0$ denotes the component meeting
$O$. The components $D_i$ have multiplicity two. The freedom of
renumbering components is killed by the following elementary
observation:

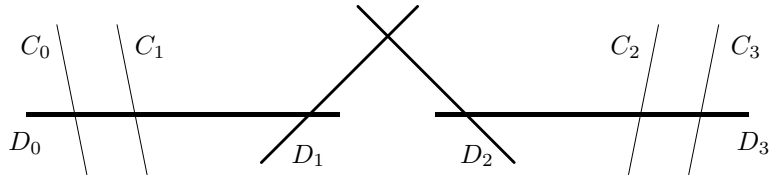
\begin{figure}[!ht]
\begin{center}
\setlength{\unitlength}{0.08mm}
\begin{picture}(1300,320)

\put(1150,0){\line(1,5){50}}
\put(1220,200){$C_3$}
\put(1050,0){\line(1,5){50}}
\put(1020,200){$C_2$}
\put(150,0){\line(-1,5){50}}
\put(250,0){\line(-1,5){50}}
\put(230,200){$C_1$}
\put(40,200){$C_0$}

{\Thicklines
\put(50,100){\line(1,0){520}}
\put(440,20){\line(1,1){260}}
\put(730,100){\line(1,0){520}}
\put(860,20){\line(-1,1){260}}}
\put(1230,40){$D_3$}
\put(20,40){$D_0$}
\put(490,20){$D_1$}
\put(770,20){$D_2$}
\end{picture}
\end{center}
\caption{The fibre of type $I_3^*$ at 0} \label{Fig-3}
\end{figure}

\begin{Lemma}\label{Lem:P}
Up to renumbering, the section $P$ meets the components $C_3$ and
$\Theta_9$.
\end{Lemma}

Lemma \ref{Lem:P} enables us to describe $P$ in terms of the
generators of the trivial lattice $V$. Since $P$ is given by
polynomials (of low degree), it does not meet $O$. We will use the
$\Q$-divisors
\begin{eqnarray*}
A & = & \frac{1}{4}(\Theta_1+2\Theta_2+\hdots+9\Theta_9+6\Theta_{10}+3\Theta_{11})\;\;\text{ and}\\
B & = & \frac{1}{4}(2C_1+4D_0+6D_1+8D_2+10D_3+5C_2+7C_3).
\end{eqnarray*}

\begin{Corollary}\label{Cor:P}
In $V\otimes_{\Z}\Q$, the section $P$ is given as
\[
P=O+2F-A-B.
\]
\end{Corollary}

\emph{Proof:} Since $4\,P\equiv 0\mod V$, it suffices to check the following intersection numbers:
\[(P.O)=0, (P.F)=1, (P.\Theta_i)=\delta_{i,9}, (P.D_j)=0, (P.C_j)=\delta_{j,3} \text{ and } (P^2)=-2.\]

\section[$X/{\F}_2$]{$\boldsymbol{X/{\F}_2}$}
\label{s:X/2}

In this section, we shall consider the elliptic K3 surface
$X/\F_2$. This coincides with the reduction $X_2$ of $X$ at 2.
Reducing equation (\ref{eq:K3}) modulo 2, we obtain the affine
equation
\begin{eqnarray}\label{eq:X/2}
X/\F_2:\;\;\;y^2+s^2\,x\,y=x^3+s^2\,x.
\end{eqnarray}
Recall that this has only two singular fibres. They sit at 0 and
$\infty$ and have types $I_3^*$ and $I_{12}$. Since 2 stays inert
in $\Q(\sqrt{-3})$, the local $L$-factor is given by
\[
P_2(X/\F_2,T)=(1-2\,T)^{21}\,(1+2\,T)
\]
because of Corollary \ref{Cor:zeta}. In accordance with the Tate
conjecture, we shall prove that $X$ is a supersingular K3 surface.
To be precise, we claim
\begin{Proposition}\label{Prop:Picard-X_2}
$X/\F_2$ is a supersingular K3 surface with Picard numbers
\[
\rho(X/\F_2)=21 \;\;\text{ and }\;\; \rho(X/\F_4)=22.
\]
\end{Proposition}
This proposition completes the proof of Proposition
\ref{Prop:Tate}. It will be established by finding explicit
generators for $NS(X/\overline\F_2)$. Additionally to the
reduction of $NS(X/\Q)$, we need two generators. By the formula of
Shioda-Tate, these can be given as sections of the elliptic
fibration $X/\overline\F_2$.

In detail, we compute some additional sections of
$X/\overline\F_2$ which are not derived from $MW(X/\Q)$ by way of
reduction. Then we determine two among them which supplement
$NS(X/\Q)$ to generate $NS(X/\overline\F_2)$. Let $\alpha$ be a
generator of $\F_4$ over $\F_2$, i.e.~$\alpha^2+\alpha+1=0$. Among
others, we found the following sections. Here, the inverse refers
to the group law on the generic fibre.
\begin{center}
\begin{tabular}{rl}
 section & inverse\\
\hline
&\\
$Q=(1,1)$ & $(1,1+s^2)$\\
$(s^2,s^2)$ & $(s^2,s^2+s^4)$\\
$(s+s^3,s^3+s^4)$ & $(s+s^3,s^4+s^5)$\\
$R=(s+\alpha s^3,\alpha^2s^4+\alpha s^5)$ & $(s+\alpha s^3,s^3+\alpha^2s^4)$\\
$(1+s^4,1+\alpha s^2+\alpha^2s^6)$ & $(1+s^4, 1+\alpha^2s^2+\alpha s^6)$
\end{tabular}
\end{center}
We shall now prove that $MW(X/\overline\F_2)$ has rank 2 and can
be generated by the sections $Q$ and $R$ together with the
torsion-section $P$. This will be achieved with the help of the
\emph{height pairing} on the Mordell-Weil group, as introduced by
Shioda in \cite{ShMW}. Let $V$ denote the trivial lattice of
$X/\F_2$. This is exactly the reduction of the trivial lattice of
$X/\Q$. The height pairing is defined via the orthogonal
projection
\[
\varphi: MW(X/\overline\F_2) \to V^\bot\otimes_\Z\Q\subset
NS(X/\overline\F_2)\otimes_\Z\Q.
\]
This projection only uses the information which (simple) component
of a reducible fibre a section meets. Let $(\cdot. \cdot)$ denote
the intersection form on $NS(X/\overline\F_2)\otimes_{\Z}\Q$.
Shioda's height pairing is defined by
\begin{eqnarray*}
<\cdot,\cdot>:\;\; MW(X/\overline\F_2)\times MW(X/\overline\F_2) & \to & \;\;\;\;\;\;\;\;\;\;\;\Q\\
\mathcal P\;\;\;\;\;\;\;\;\;\;\;\;\;\;\;\;\;\;\;\; \mathcal
Q\;\;\;\;\; & \mapsto & -(\varphi(\mathcal P). \varphi(\mathcal
Q)).
\end{eqnarray*}
The height pairing is symmetric and bilinear. It induces the
structure of a positive definite lattice on
$MW(X/\overline\F_2)/MW(X/\overline\F_2)_\text{tor}$ (the
Mordell-Weil lattice). For $\mathcal P, \mathcal Q\in
MW(X/\overline\F_2)$, Shioda shows that
\[
<\mathcal P,\mathcal Q>=\chi(\OO_X)-(\mathcal P.\mathcal
Q)+(\mathcal P.O)+(\mathcal Q.O)-\sum_v \text{contr}_v(\mathcal
P,\mathcal Q).
\]
Here the sum runs over the cusps and contr$_v$ can be computed
strictly in terms of the components of the singular fibres which
$\mathcal P$ and $\mathcal Q$ meet \cite[Thm.~8.6]{ShMW} (cf.~also
\cite[(8.16)]{ShMW}). The following lemma is easily checked:

\begin{Lemma}\label{Lem:Q,R}
$Q$ meets $C_0$ and $\Theta_8$ while $R$ intersects $C_3$ and $\Theta_1$.
\end{Lemma}

We are now in the position to compute the projections $\varphi(Q)$
and $\varphi(R)$. However, we will postpone their explicit
calculation, since we will only need this for the explicit
verification of Shioda's Conjecture \ref{Conj:Shioda}. Here, we
can use Shioda's results from \cite{ShMW} to prove that
$NS(X/\bar\F_3)$ is generated by the trivial lattice $V$ and the
sections $P, Q$ and $R$.

\begin{Lemma}\label{Lem:<Q,R>}
Let $Q, R$ be the sections of $X/\F_4$ as specified above. Then
\[
\begin{pmatrix}
<Q,Q> & <Q,R>\\ <Q,R> & <R,R>
  \end{pmatrix}=\frac{1}{3}
\begin{pmatrix}
4 & 2\\2 & 4
\end{pmatrix}.
\]
\end{Lemma}

\emph{Proof:} Since they are given by polynomials (of low degree),
$Q$ and $R$ do not meet $O$. For the self intersection numbers, we
thus only miss the contributions from the singular fibres. These
are derived from \cite[(8.16)]{ShMW} with the help of Lemma
\ref{Lem:Q,R}:
\begin{eqnarray*}
<Q,Q> & = & 4-\frac{4\cdot 8}{12}=\frac{4}{3}\\
<R,R> & = & 4-\frac{7}{4}-\frac{11}{12}=\frac{4}{3}.
\end{eqnarray*}
For the remaining entry, we furthermore have to analyze the
intersection of $Q$ and $R$. We have to find the common zeroes of
\begin{eqnarray*}
\alpha s^3+s+1 & = & \alpha(s+\alpha^2)(s^2+\alpha^2s+1) \;\;\;\;\text{and}\\
\alpha s^5+\alpha^2s^4+1 & = & \alpha
(s+1)(s+\alpha^2)(s^3+\alpha^2s+1).
\end{eqnarray*}
The factors in the above factorization are irreducible over
$\F_4$. Hence, the only common zero is $s=\alpha^2$. Since this
occurs with multiplicity one, the intersection is transversal.
Hence, we deduce
\[
<Q,R> = 2-1-\frac{1\cdot 4}{12}=\frac{2}{3}.
\]
This finishes the proof of Lemma \ref{Lem:<Q,R>}.

\begin{Proposition}\label{Prop:X_2}
$NS(X/\overline\F_2)$ has discriminant $-4$. It is generated by
the trivial lattice $V$ and the sections $P, Q$ and $R$, so
$MW(X/\overline\F_2)=MW(X/\F_4)=<P,Q,R>$.
\end{Proposition}
\emph{Proof:} Consider the lattice $N$ which is generated by the
trivial lattice $V$ and the sections $P, Q$ and $R$. Clearly, this
is a sublattice of $NS(X/\overline\F_2)$. We can identify
$N'=<V,P>\cong NS(X/\Q)\subset N$. Recall that this sublattice has
discriminant $-3$.

Then we use the orthogonal projection $\varphi$ from
$NS(X/\overline\F_2)$ to $V^\vee\otimes_{\Z}\Q$. Since $P\in
V\otimes_{\Z}\Q$, it follows that
\[
\text{discr}\, N = (\text{discr} \,N')
\;\,\text{det}\begin{pmatrix}
(\varphi(Q). \varphi(Q)) & (\varphi(Q). \varphi(R))\\
(\varphi(Q). \varphi(R)) & (\varphi(R). \varphi(R))
  \end{pmatrix}.
  \]
By Lemma \ref{Lem:<Q,R>}, the matrix has determinant $\frac 43$.
Hence, $N$ has discriminant $-4$ and in particular rank 22. Since
both values are maximal possible (cf.~Thm.~\ref{Thm:Artin-inv}),
we obtain $N=NS(X/\overline\F_2)$. The claim
$MW(X/\overline\F_2)=<P, Q, R>$ then follows from the formula of
Shioda-Tate. This proves Proposition \ref{Prop:X_2}.

Proposition \ref{Prop:X_2} implies Proposition
\ref{Prop:Picard-X_2}. It verifies the Tate conjecture for
$X/\F_{2^r}$ with any $r\in\N$. We also deduce the validity of
Theorem \ref{Thm:Artin}. Furthermore, we have seen that $X_2$
has Artin invariant $\sigma_0=1$.

We conclude this section by verifying Conjecture \ref{Conj:Shioda}
for $X/\Q$ and its reduction at 2. Note that it suffices to
consider $NS(X/\Q)$ and $NS(X/\F_4)$. The dual of $NS(X/\Q)$ will
always refer to the embedding into $NS(X/\F_4)$.

To give the explicit form of $\varphi(Q)$ and $\varphi(R)$, recall
the $\Q$-divisors $A$ and $B$ from Section \ref{s:NS}. The
projections are easily computed as
\begin{eqnarray*}
\varphi(Q) & = & Q-O-2F+\frac{1}{3}(\Theta_1+2\Theta_2+\hdots+8\Theta_8+6\Theta_9+4\Theta_{10}+2\Theta_{11})\\
\varphi(R) & = &
R-O-2F+B+\frac{1}{12}(11\Theta_1+10\Theta_2+\hdots+\Theta_{11}).
\end{eqnarray*}
Since the denominators are divisible by $3$, $\varphi(Q),
\varphi(R)\not\in NS(X/\F_4)$. In other words,
\begin{eqnarray}\label{eq:neqq}
<\varphi(Q),\varphi(R)>\, \supsetneqq NS(X/\Q)^\bot\,\supseteq
<3\varphi(Q),3\varphi(R)>.
\end{eqnarray}
We claim that
\begin{eqnarray}\label{eq:NS^}
NS(X/\Q)^\bot=<3\varphi(Q),
\varphi(Q)+\varphi(R)>=<3\varphi(R),\varphi(Q)+\varphi(R)>.
\end{eqnarray}
Since the ranks are two, the inequality of (\ref{eq:neqq}) implies
that the claim is equivalent to
\[
\varphi(Q)+\varphi(R)\in NS(X/\F_4).
\]
Since $P\equiv -A-B \mod V$, it suffices to show that
\[
\varphi(Q)+\varphi(R)\equiv A+B\mod{<V,Q,R>}.
\]
Explicitly, we have
\begin{eqnarray*}
\varphi(Q)+\varphi(R) & = & Q+R-2O-4F+\frac{1}{12}(11\Theta_1+10\Theta_2+\hdots+\Theta_{11})
\\
&& +B+\frac{1}{3}(\Theta_1+2\Theta_2+\hdots+8\Theta_8+6\Theta_9+4\Theta_{10}+2\Theta_{11})\\
& \equiv & B +\frac{1}{12}\sum_j(12-j)\Theta_j+\frac{1}{3}\sum_jj\Theta_j\\
& = & B+\frac{1}{12}\sum_j(12-j+4j)\Theta_j\\
& \equiv & B+\frac{1}{4}\sum_jj\Theta_j \equiv A+B \;\;\mod{<V,Q,R>}.
\end{eqnarray*}
This proves claim (\ref{eq:NS^}). We shall now verify Conjecture
\ref{Conj:Shioda} for $X/\Q$ and its reduction at 2. Here
$T_2=<3\varphi(Q), \varphi(Q)+\varphi(R)>$ with intersection form
\[
\begin{pmatrix}
(3\varphi(Q). 3\varphi(Q)) & (3\varphi(Q). \varphi(Q)+\varphi(R))\\
(3\varphi(Q). \varphi(Q)+\varphi(R)) & (\varphi(Q)+\varphi(R).
\varphi(Q)+\varphi(R))
  \end{pmatrix}\hspace{1cm}\]\[\hspace{1cm}=-
\begin{pmatrix}
12 & 6\\6 & 4
\end{pmatrix}=-2
\begin{pmatrix}
6 & 3\\3 & 2
\end{pmatrix}\sim -2
\begin{pmatrix}
2 & 1\\1 & 2
\end{pmatrix}.
\]
Hence, Conjecture \ref{Conj:Shioda} holds with $m=2$.

\begin{Remark}
In \cite{DK}, Dolgachev and Kond\=o give several models of the supersingular K3 surface in characteristic $2$ with Artin invariant $\sigma_0=1$. The most natural might be the quasi-elliptic fibration
\[
Z:\;\;y^2 = x^3 + t^2 x +t^{11}.
\]
The only exceptional fibre of this fibration has type $I_{16}^*$. This shows that
\[
NS(Z) = U \oplus D_{20}.
\]
We realize a $\bar\F_2$-isomorphic fibration with the above exceptional fibre in terms of our model $X/\F_2$ by giving an effective divisor $L$ of type $I_{16}^*$. Denote by $\iota R$ the conjugate of the section $R$ in $\F_4$ and by $[3P]$ the section $(-s,s^3)$ (the inverse of $P$ with respect to the group law on the fibre). Recall that both, $R$ and $\iota R$ meet $C_3$ and $\Theta_1$. $[3P]$ meets $C_2$ and $\Theta_3$ (cf.~Fig.~\ref{fig:3IV*}). Thus the required divisor $L$ (over $\F_2$) is given by
\begin{eqnarray*}
L & = & R+\iota R
+2(\Theta_1+\Theta_0+\Theta_{11}+\Theta_{10}+\Theta_9+\Theta_8+\Theta_7+\Theta_6+\Theta_5\\
&& \;\;\; +\Theta_4+\Theta_3+[3P]+C_2+D_3+D_2+D_1+D_0)
+C_0+C_1.
\end{eqnarray*}

\end{Remark}

\section{Twisting}
\label{s:twist-X}

This section concludes the investigation of the extremal elliptic
K3 surface $X$ by commenting on twisting. Our aim is to show that
the twists of the associated newform $f$ with rational
coefficients are in correspondence with twists of $X$. This is
non-trivial, since $f$ admits cubic twists.

For the quadratic twisting, this is well-known. For instance, it
follows from point counting that the twist of $X$ over
$\Q(\sqrt{d})$ corresponds to the quadratic twist
$f\otimes\left(\frac{d}{\cdot}\right)$. Note that there is a model
of this twist over $\Q$ with bad reduction exactly at 3 and the
primes dividing the discriminant of $\Q(\sqrt{d})$.

We shall now come to the cubic twists. For $X$, they can be
achieved by considering the base change
\[
\pi_d: s \mapsto d\, s^3
\]
instead of the original $\pi$. Here, we can restrict to positive
cube-free $d$. Let $X^{(d)}$ denote the pull-back from the
rational elliptic surface $Y$ by $\pi_d$. An affine model is given
by
\begin{eqnarray}\label{eq:X^{(d)}}
X^{(d)}:\;\; y^2+d^2s^2\,x\,y=x^3+2s\,x^2+s^2\,x.
\end{eqnarray}
Of  course, $X^{(d)}$ has the same configuration of singular
fibres as $X$ and also inherits the trivial action of
Gal$(\overline\Q/\Q)$ on $NS(X^{(d)})$. By construction, it has
bad reduction exactly at the prime divisors of $3d$.

For instance, consider the twist $X^{(3)}$. Since this has good
reduction away from $3$, the arguments of Section \ref{s:new}
apply to determine the associated newform. Recall the
Gr\"o\ss{}encharakter $\psi_f$ associated to the newform $f$.
Counting points mod $7$, we find that $X^{(3)}$ corresponds to the
twist $\psi_f\otimes\left(\frac{3}{\cdot}\right)_3$ (with trace
$11$ at $7$). We now give the general statement:

\begin{Theorem}\label{Thm:X^{(d)}}
$X^{(d)}$ corresponds to the twist
$\psi_f\otimes\left(\frac{d}{\cdot}\right)_3$, i.e.
$L(T_{X^{(d)}},s)=L(\psi_f\otimes\left(\frac d\cdot\right)_3,s).$
\end{Theorem}

The proof of this theorem is organized as follows: We exhibit
another extremal elliptic fibration on the surfaces $X^{(d)}$.
This will have fibres with CM by $\Q(\sqrt{-3})$ such that we can
perform the twisting fibrewise. Then the claim will follow.

\begin{figure}[!ht]
\begin{center}
\setlength{\unitlength}{0.1mm}
\begin{picture}(1500,1000)

\imagerotate{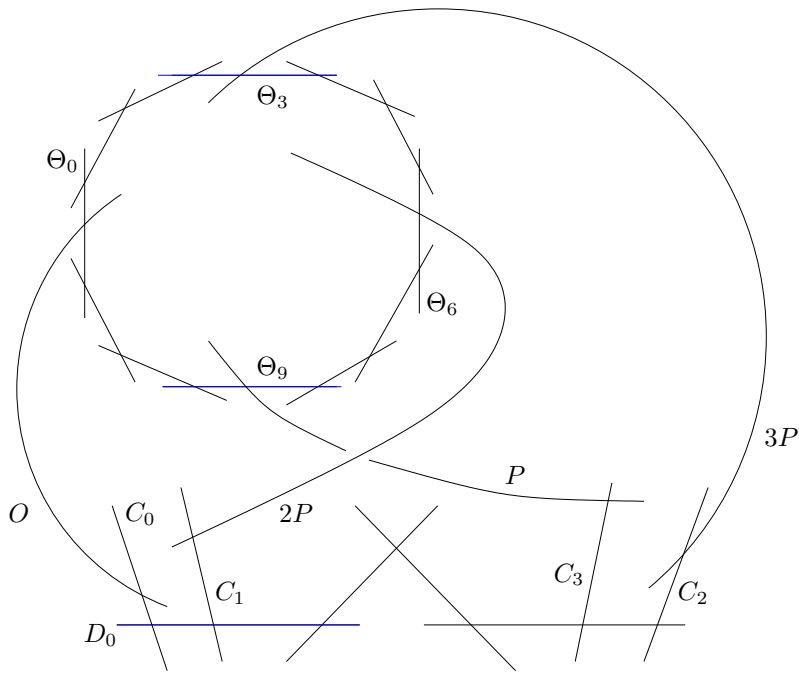}{150,0}{100mm}{3IV}{150,0}{100mm}
\put(140,200){$O$} \put(500,200){$2P$} \put(800,250){$P$}
\put(1145,300){$3P$} \put(240,40){$D_0$} \put(295,200){$C_0$}
\put(415,95){$C_1$} \put(865,120){$C_3$} \put(1030,95){$C_2$}
\put(190,670){$\Theta_0$} \put(695,480){$\Theta_6$}
\put(470,395){$\Theta_9$} \put(470,755){$\Theta_3$}
%\put(150,0){\epsfig{file=3IV*.eps,width=100mm}}
%\put(200,0){\epsfig{file=dessin1,width=20mm}}
%\put(800,40){\epsfig{file=dessin2,width=28mm}}
%\put(100,-70){\epsfig{file=versuch1,width=42mm}}
%\put(750,-50){\epsfig{file=versuch2,width=40mm}}

\end{picture}
\end{center}
\caption{The fibration on $X^{(d)}$ with three fibres of type
$IV^*$} \label{fig:3IV*}
\end{figure}

We want to give an elliptic fibration on $X^{(d)}$ with the
component $D_0$ as a section. Therefore, we work with the affine
equation of the first blow-up of the Weierstrass model
(\ref{eq:X^{(d)}}) at the $D_3$-singularity:
\begin{eqnarray}\label{eq:fibr-IV}
y'^2+d^2\,s^2x'y'=sx'(x'+1)^2.
\end{eqnarray}
Then $D_0=\{s=y'=0\}$, so the fibration is the affine projection
on the $x'$-coordinate. We employ the usual notation, i.e.~replace
$x'$ by $t$ and likewise for $s$ and $y'$. Then equation
(\ref{eq:fibr-IV}) becomes
\begin{eqnarray*}%\label{eq:3IV*}
y^2+d^2\,tx^2y=t(t+1)^2x.
\end{eqnarray*}
In order to obtain a projective model for this, we first homogenize
\[
s^3y^2z+d^2s^2tx^2y=t(t+s)^2xz^2.
\]
Then the change of variable $y\mapsto \frac{t+s}{s}y$ gives
\begin{eqnarray}\label{eq:homog}
X^{(d)}:\;\;s(t+s)y^2z+d^2stx^2y=t(t+s)xz^2.
\end{eqnarray}
This has six singularities, two in each fibre above $0, -1$ and
$\infty$. Their resolution produces three fibres of type $IV^*$.
These can be read off directly from the original fibration. We
sketch this in Figure \ref{fig:3IV*}.

\vspace{0.2cm}

The sections of the new fibration are $D_0, \Theta_3$ and
$\Theta_9$, as printed in blue, and the remaining components form
the three singular fibres of type $IV^*$. In particular, they are
all defined over $\Q$, as we have seen before.

The smooth fibres of the above fibration are elliptic curves with
CM by $\Q(\sqrt{-3})$. However, the impact of twisting on the
associated newforms is not visible so far. Therefore, we transform
eq.~(\ref{eq:homog}) into Weierstrass form. The procedure from
\cite[\S 8]{Ca} gives
\begin{eqnarray*}
X^{(d)}:\;\;y^2 + d^2 s^2t^2(s+t)^2 y = x^3.
\end{eqnarray*}
By \cite[18, Thm.~4]{IR}, the Gr\"o\ss{}encharaktere associated to
the smooth fibres are twisted by $\left(\frac {d}{\cdot}\right)_3$
upon moving from $X^{(1)}$ to $X^{(d)}$. Using the Lefschetz fixed
point formula, we can express the trace of Frobenius on
$T_{X^{(d)}}$ as the sum of the traces on the smooth fibres. Hence
Theorem \ref{Thm:X^{(d)}} follows.

\begin{Remark}
From the classification of \cite{S-CM}, it follows that any
newform of weight 3 with rational coefficients and CM by
$\Q(\sqrt{-3})$ can be realized geometrically by some twist of
$X$.

For the only other CM-field with non-quadratic twists,
$\Q(\sqrt{-1})$, the corresponding statement can be established
using the Fermat quartic in $\PP^3$.
\end{Remark}

\begin{Corollary}
Up to the bad Euler factors, we have
\[
\zeta(X^{(d)}, s) \circeq \zeta(s)\; \zeta(s-1)^{20}\;
L(\psi_f\otimes\left(\frac d\cdot\right)_3,s)\; \zeta(s-2).
\]
\end{Corollary}

\vspace{0.8cm}

{\bf Acknowledgement:} This paper is condensed from the last
chapter of my PhD-thesis. I would like to express my gratitude to
my advisor K.~Hulek for all his support and encouragement during
this period. I thank B.~van Geemen and T.~Shioda for many useful
comments and discussions.\\
Generous support from the DFG-Schwerpunkt 1094 "Globale Methoden
in der komplexen Geometrie" and the network Arithmetic Algebraic
Geometry, a Marie Curie Research Training Network, is gratefully
acknowledged. The paper was finished while I enjoyed the
hospitality of the Dipartimento di Matematica "Frederico Enriques"
of Milano University. I am very grateful for the stimulating
atmosphere.

\vspace{0.3cm}

Matthias Sch\"utt\\ 
Department of Mathematical Sciences\\
University of Copenhagen\\
Universitetsparken 5\\
2100 Copenhagen\\
Denmark\\
mschuett@math.ku.dk

\end{document}